\documentclass[11pt]{article}
\usepackage{amsmath,amssymb,amsbsy,amsthm,fullpage}

\usepackage[utf8]{inputenc} % allow utf-8 input
\usepackage[T1]{fontenc}    % use 8-bit T1 fonts
\usepackage{hyperref}       % hyperlinks

\let\epsilon\varepsilon
\let\phi\varphi

\def\X{{\cal X}}
\def\x{{\bf  x}}
\def\N{\mathbb N}

\def\C{\mathcal C}
\def\E{{\bf E}}

\def\-as{\text{-a.s.}}

\def\bd{{ D}}

\newtheorem{theorem}{Theorem}\newtheorem*{theorem*}{Theorem}
\newtheorem{definition}{Definition}

\newtheorem{proposition}{Proposition}
\newtheorem{corollary}{Corollary}

\begin{document}
\title{Universality of Bayesian mixture predictors}
\author{Daniil Ryabko}
% \author{{\bf Daniil Ryabko}\\    daniil@ryabko.net, \\  INRIA 
% }
\date{}
%\date{draft paper, do not distribute}
\maketitle

\begin{abstract}
The problem is that of sequential probability forecasting for finite-valued time series. The data is generated by 
an unknown  probability distribution over the space of all one-way infinite sequences.
It is known that this measure belongs to a given set $\C$, but the latter is completely arbitrary 
(uncountably infinite, without any structure given).  The performance is measured 
with  asymptotic average log loss.  In this work it is shown that the minimax asymptotic
performance  is always attainable, and it is attained by a convex combination of a countably many 
measures from the set $\C$ (a Bayesian mixture). This was previously only known for the case when the best achievable asymptotic error is 0.
This also contrasts previous results that show that in the non-realizable case all Bayesian mixtures may be suboptimal, while 
there is a predictor that achieves the optimal performance.
\end{abstract}

% \begin{keywords}
%   Sequence Prediction,  Time Series, Online Prediction, Bayesian Prediction
% \end{keywords}
\section{Introduction}

Given a  sequence $x_1,\dots,x_n$ of observations $x_i\in\X$, where $\X$ is a finite set, we 
want to predict what are the probabilities of observing $x_{n+1}=x$ for each $x\in\X$, before $x_{n+1}$ is revealed,
after which the process continues sequentially.
The sequence $x_1,\dots,x_n,\dots$ is generated by some unknown stochastic process $\mu$, a probability measure
on the space of one-way infinite sequences $\X^\infty$. 
Further, a set $\C$ of such measures is given, and it is known that $\mu\in\C$.  The set $\C$ can be thought of as the set of models, experts, or the set of strategies of the adversary (a.k.a.\ Nature).
The requirement that the true measure $\mu$ that generates the data is in $\C$ means that we are in the ``realizable'' case 
of the problem (in other words, there is  at least one expert that is optimal all the time).

Unlike most of the literature on the subject, which assumes that the set $\C$ is parametrized and endowed with some structure,
here we would like to treat the problem in full generality and thus shall not make any assumptions whatsoever on the set $\C$
or its elements.  Note that making even such innocuous-looking assumptions on the set $\C$ as are the common topological ones, 
 such as local compactness, separability, tightness, not to mention stronger assumptions involving the existence of densities or smoothness, implicitly gives the problem a structure (e.g., a topology in which the assumption is fulfilled) that in itself constitutes a large part of the solution. 
Here we are interested rather in the fundamental question of which principles to use when  choosing  a model 
for a problem, and thus would like not to make any assumptions at all (not even measurability). It is also worth reiterating that the measures in $\C$ are
not required to be i.i.d., finite-memory, mixing, etc.

We are interested in the question whether it is possible to attain the minimax optimal asymptotic performance   by using a combination 
of measures in $\C$ as a predictor. A combination is a measure  of the form $\nu=\int_\C dW$ where $W$ is some measure over $\C$ 
(or, more generally, over a measurable subset of $\C$; recall that $\C$ itself is not required to be measurable). The measure 
$W$ can be thought of as a prior distribution, and prediction is then by evaluating the posterior $\nu(\cdot|x_1,\dots,x_n)$ on 
the observed sequence $x_1,\dots,x_n$.  In other words,  we are asking whether it is possible to achieve optimal prediction with a Bayesian predictor with some prior (whether there exists such a prior); however, we are not interested in optimality with prior probability 1, but rather in the  minimax (worst-case) asymptotic optimality of such predictors. The answer we obtain is positive. 

Thus, {\bf the main result} is that   the minimax asymptotic performance is always attainable and it is attained by a  combination of  countably many measures from the set $\C$.  Note that this statement  is completely assumption-free: not even measurability of $\C$ is required.
% : it shows that the best possible
% performance in a prediction problem is {\em always} attained by a Bayesian mixture predictor. 

Previously, this result was only established, in \cite{Ryabko:10pq3+}, under  the assumption that
 there is a predictor whose error is asymptotically 0 on any measure $\mu\in\C$, that is, the minimax asymptotic error is 0. 
 Here we get rid of this (last) assumption.  Examples of cases where the minimax asymptotic error is greater than 0 are easiest
to come by if we suppose that some aspects of the process are completely arbitrary. The easiest example is when nothing is predictable: the data
is an arbitrary deterministic sequence. This example gives the maximal possible worst-case asymptotic error of $\log|\X|$. 
A more meaningful example is that of processes with (frequent) abrupt changes:
between the changes, the distribution belongs to some (nice) given family (e.g., i.i.d.\ Bernoulli trials) but when a changes occur is not known, and a change is  to an arbitrary distribution in the family. This example is considered in more detail in Section~\ref{s:ex}.

Moreover, the case when the best possible asymptotic error is greater than 0 is particularly important in light of 
recent results achieved in the non-realizable case,
% Let us compare this result to what happens in the non-realizable case, 
that is, when the measure $\mu$ generating 
the data does not have to belong to $\C$. In this case, one is interested in the regret with respect to $\C$, that is, 
the performance of the predictor minus the best performance of all the measures in $\C$ on the given $\mu$. 
It  has recently been shown in \cite{Ryabko:16pqnot}  that, in the non-realizable case, 
it can happen that the best regret a predictor can achieve with respect to a set $\C$ is zero, but 
any Bayesian mixture predictor has regret  bounded away from 0 by a large constant (see Section~\ref{s:not} for a precise formulation).
In other words, the experts in $\C$ are useless: one can do as well as any of them, but not by combining them. 
 Note that any such set $\C$ has to 
be uncountable, which brings it out of the traditional expert advice settings (a survey of which can be found in \cite{Cesa:06}).

Putting these results together, we reach the following  fundamental {\bf recommendation for choosing a model} for sequential data:

{\em Better take a model large enough to make sure it includes the process that generates the data, even if it makes
 the worst-case asymptotic error  larger than zero, for otherwise any combination of predictors in the model 
class may be useless.}

To the initiated reader this  result has a distinct decision-theoretic flavour to it. Indeed, as we explain in Section~\ref{s:dec},
it provides (a strong form) of the the complete-class theorem for the problem of sequential prediction, as well as a 
partial form of the minimax theorem.

\noindent{\bf Related work.}
The literature on (nonparametric) sequential prediction is huge, and we do not attempt to provide an adequate survey here.
Some pointers to climb references from in different branches of science are: \cite{Cesa:06} for the expert advice setting (machine learning side),
\cite{Kalai:94,Noguchi:15} for the non-parametric Bayesian approach (econometrics side; most results with prior probability~1), 
\cite{BRyabko:88,Morvai:97,Gyorfi:07} for predicting stationary ergodic time series (perhaps the largest class considered in statistics; non-parametric statistics/information theory side), \cite{Solomonoff:78,Hutter:04uaibook} for predicting computable measures. 
%, as well as \cite{Ryabko:10pq3+,Ryabko:11pq4+} for a general assumptionless approach attempting to unify the above.
%
 The study of the realizable and non-realizable sides of the prediction problem together in the setting considered here  has been initiated in \cite{Ryabko:11pq4+} that    also poses  the question that \cite{Ryabko:16pqnot}  resolves. 

To author's knowledge, this is the first work to consider the general case case when  the best achievable asymptotic regret is greater than 0. One specific 
example that was  considered before  is that of processes with abrupt changes mentioned above. The work \cite{Willems:96} considers the case 
when the processes between the changes are i.i.d., but the method proposed is general. It is, in fact, a Bayesian construction, where the prior is 
over all possible sequences of changes. The goal is to minimize the regret  with respect to the predictor that knows the sequence of changes (but not the distributions); the best achievable asymptotic regret was not considered directly. Subsequent work on this problem is largely devoted to computational considerations;
see \cite{Gyorgy:12} and references therein. 

Related  decision-theoretic results concern the setting of the problem for ``predicting'' just one (the first) 
 symbol of the sequence. % (but note that in the minimax setting there is not much difference between just one and any finitely many symbols), Thus, the minimax theorem for the settings where instead of asymptotic 
% average loss we are concerned with just the KL divergence or another measure of discrepancy 
%between the distributions over one step.
 For KL divergence (expected log loss) these results include \cite{BRyabko:79,Gallager:76,Haussler:97}; 
 a variety of generalizations to other losses is presented in \cite{Grunwald:04}.

% The rest of the paper is organized as follows. The next section introduces notation and definitions. 
% Section~\ref{s:main} presents the main result. Section 

\section{Preliminaries}\label{s:pre}
Let $\X$ be a finite set. The notation $x_{1..n}$ is used for $x_1,\dots,x_n$. 
 We consider (probability) measures on $(\X^\infty,\mathcal F)$, where $\mathcal F$
is the usual Borel sigma-field generated. 
%  by the cylinder sets  $[x_{1..n}]$, $x_i\in\X, n\in\N$, 
% where $[x_{1..n}]$ is the set of all infinite sequences that start with $x_{1..n}$. 
% Since we are only interested in those measures on  $(\X^\infty,\mathcal F)$ which are probability measures (the measure of $\X^\infty$ equals 1),
% we call them simply {\em measures}.
For a  finite set $A$ denote $|A|$ its cardinality.
We use  $\E_\mu$ for expectation with respect to a measure $\mu$.

For two measures $\mu$ and $\rho$ introduce the {\em expected cumulative Kullback-Leibler divergence (KL divergence)} as
\begin{equation}\label{eq:akl} 
  d_n(\mu,\rho):=  \E_\mu
  \sum_{t=1}^n  \sum_{a\in\X} \mu(x_{t}=a|x_{1..t-1}) \log \frac{\mu(x_{t}=a|x_{1..t-1})}{\rho(x_{t}=a|x_{1..t-1})}=
  \sum_{x_{1..n}\in\X^n}\mu(x_{1..n}) \log \frac{\mu(x_{1..n})}{\rho(x_{1..n})}.
\end{equation}
In words, we take the expected (over data) average (over time) KL divergence between $\mu$- and $\rho$-conditional (on the past data) 
probability distributions of the next outcome; and this gives simply the $\mu$-expected log-ratio of the likelihoods.

\begin{definition}
Define the asymptotic average KL loss of $\rho$ on $\mu$ as 
\begin{equation}\label{eq:dn}
\bd(\mu,\rho):=\limsup_{n\to\infty}{1\over n}d_n(\mu,\rho).
\end{equation}
For a set $\C$ of measures define 
$$
\bd(\C,\rho):=\sup_{\mu\in\C}\bd(\mu,\rho).
$$
\end{definition}

% A simple but useful identity that we will need 
% is the following
% \begin{equation}\label{eq:kl}
%  d_n(\mu,\rho)=-\sum_{x_{1..n}\in\X^n}\mu(x_{1..n}) \log \frac{\rho(x_{1..n})}{\mu(x_{1..n})},
% \end{equation}
% where on the right-hand side we have simply the KL divergence between measures $\mu$ and $\rho$ restricted to the first $n$ observations.
The main quantity of interest is the following minimax loss:
\begin{definition}
For a set $\C$ of measures define 
\begin{equation}\label{eq:vc}
  V_\C:=\inf_{\rho\in\mathcal P} \bd(\C,\rho)=\inf_{\rho\in\mathcal P}\sup_{\mu\in\C}\bd(\mu,\rho),
\end{equation}
where the infimum is taken over the set $\mathcal P$ of all probability measures on $(\X^\infty,\mathcal F)$.
\end{definition}
Thus, $V_\C$ is the minimax loss for the set $\C$ of strategies of the Nature and unrestricted set of statistician's strategies. 

\section{Main result}\label{s:ba}

The main result shows that the minimax loss is always achievable and is achieved by a convex combination of measures in $\C$~---
without any assumptions on $\C$. Moreover, for any predictor $\rho$ there is a convex combination of measures in $\C$ that is as good as $\rho$.

\begin{theorem}\label{th:2} For any set  $\C$ of probability measures on $(\X^\infty,\mathcal F)$,
there exist a sequence of measures $\mu_k\in\C$, $k\in\N$ and a sequence of real weights  $v_k>0,k\in\N$ whose sum is 1, such that
for the measure $\phi:=\sum_{k\in\N} v_k\mu_k$ we have 
$$
\bd(\C,\phi)= V_\C.
$$
Moreover, for every measure $\rho$ there exists a predictor $\phi$ of the form above such that $D(\mu,\phi)\le D(\mu,\rho)$ for all $\mu\in\C$.
\end{theorem}
Before giving the proof, we present   informally some ideas behind it. 
Imagine first that for the set $\C$ we already knew  a predictor $\rho$ that attains the value $V_\C$.
Imagine furthermore, that for each $\mu\in\C$ the limit $\lim_{n\to\infty}{1\over n}\log\frac{\mu(x_{1..n})}{\rho(x_{1..n})}$ exists 
for $\mu$-almost all $\x=x_1,\dots,x_n,\dots\in\X^\infty$. Then we could define ($\mu$-almost everywhere) the function  $f_\mu(\x)$ whose value at $\x$ equals this limit. Let us call it the ``log-density'' function. (The reader can recognize behind the log 
the expression $\lim_{n\to\infty}\frac{\mu(x_{1..n})}{\rho(x_{1..n})}$ that  defines the familiar  densities.) Furthermore, nothing forbids us to imagine that this log-density is measurable. What we would be looking for thence is to find a countable dense subset of the set of  log-densities of all measures from $\C$. The measures $\mu$ corresponding to each log-density in this countable set would then constitute the sequence whose existence the theorem asserts. To find such a dense countable subset we could employ a standard procedure: approximate all log-densities by step functions with finitely many steps. The main technical argument is then to show that, for each level of the step functions,
 there are not too many of these functions whose steps are concentrated on different sets of  non-negligible probability,  for otherwise the requirement that $\rho$ attains $V_\C$ would be violated. Here ``not too many'' means  exponentially many with the right exponent (the one corresponding to the step of the step-function with which we approximate the density), 
and ``non-negligible probability'' means a probability bounded away (in $n$) from 0. Getting back to reality, 
we cannot say anything about the existence of the limits. What we do instead is use the step-functions approximation at each time step $n$. Since there are only countably many time steps, the result is still a countable set of measures $\mu$ from $\C$. Finally, we
are not given a measure $\rho$ that attains the value $V_\C$; so we find a sequence of predictors $\rho_{\epsilon_n}$ 
that approach this value instead, and  perform the procedure above for each $\rho_{\epsilon_n}$.

It is worth noting that the proof that \cite{Ryabko:10pq3+} obtains for the special case $V_\C=0$,
does not directly generalize. In fact, tidying up the constants in the proof in \cite{Ryabko:10pq3+}, one only obtains the asymptotic loss of $2V_\C$ for the mixture predictor presented there. It is not a problem for the case $V_\C=0$, but of course is not what we want in the 
general case.  The reason behind this problem  is that for the construction in that proof one can only use the fact that each of the measures $\mu_k$ in the sequence is as good as the predictor $\rho$ whose existence is assumed (the one that attains $V_\C=0$). 
In contrast, in the proof below we are able to use the fact that each measure in the sequence is in fact much better than $\rho$ on some subsets of $\X^n$.

% The reader could of course recognizer behind the imaginary log-densities above the familiar density functions defined 
% by  $\lim_{n\to\infty}\frac{\mu(x_{1..n})}{\rho(x_{1..n})}$, which indeed exist $\mu$-a.s.\ but is only finite
% on the part of the space $\X^\infty$ where $\mu$ is absolutely continuous with respect to $\rho$.

\begin{proof}
Define the weights $w_k:=w/k\log^2k$, where $w$ is the normalizer such that $\sum_{k\in\N}w_k=1$.
Introduce the notation $M:=\log|\X|$.

When speaking about measures $\nu$ that we construct as countable convex combinations of measures in $\C$
we will assume w.l.o.g.\  
\begin{equation}\label{eq:boundedness}
  -\log\nu(x_{1..n})\le nM+1  \text{ for all $n\in\N$
and $x_{1..n}\in\X^n$}.                      
\end{equation}
% (the value of the constant will not be important, but one can take $const=1$);
 Thus, in particular, $d_n(\mu,\nu)\le nM+1$ for all $\mu$. This boundedness can always be achieved
by taking $\nu:=1/2(\nu+p)$, where $p$ is the i.i.d.\ measure with equal probabilities of outcomes, in the case $p\in\C$,
and if $p\notin\C$ the required boundedness can be obtained as described in \cite[end of the proof (Step~{\em r}) of Theorem~5]{Ryabko:10pq3+}. The argument is technical and of no great importance so we do not reproduce it here.
%\noindent{\em Step 1} 
% For all $\mu\in\C$, $n\in\N$, $\gamma\in(0,1)$ define the sets 
% \begin{equation}\label{eq:t}
% T_{\mu,\gamma}^n:=\left\{x_{1..n}\in \X^n: \mu(x_{1..n})\ge2^{\gamma n}\rho(x_{1..n})\right\}.
% \end{equation} 

% For each $k\in\N$ let $U_k$ be a fixed finite $2^{-k}$-net in $[0,1]$ such as
%  $U_k:=\{u_k^i:i=1..2^{k}\}$ 
% where $u^i_k=[(i-1)2^{-k},i2^{-k})$. 

We start with the second statement of the theorem. Take any predictor $\rho$.
We shall  find a measure $\nu$ of the form $\sum_{k\in\N}w_k'\mu_k$, where $\mu_k\in\C$ 
such that  
\begin{equation}\label{eq:nuro}
D(\mu,\nu)\le D(\mu,\rho) \ \forall \mu\in\C.
\end{equation}

For each $\mu\in\C$, $n\in\N$ define the sets 
\begin{equation}\label{eq:tt}
T_{\mu}^n:=\left\{x_{1..n}\in \X^n:  \frac{\mu(x_{1..n})}{\rho(x_{1..n})}\ge {1\over n}\right\}.
\end{equation} 
From Markov inequality, we obtain 
\begin{equation}\label{eq:tm}
\mu(\X^n\backslash T_{\mu}^n)\le 1/n.
\end{equation}

For each $k\in\N$ let $U_k$ be the partition of  $[-\frac{\log n}{n},M+{1\over n}]$ into   $k$ intervals defined as follows.
 $U_k:=\{u_k^i:i=1..k\}$, where
$$
 u^i_k=\left\{ \begin{array}{ll}
                \left[-\frac{\log n}{n},{iM\over k}\right] & i=1, \\
 		\left(\frac{(i-1)M}{k},\frac{iM}{k}\right] & 1<i<k, \\
 		\left(\frac{(i-1)M}{k},M+\frac{1}{n}\right] & i=k.
               \end{array}
  \right.
$$
Thus, $U_k$ is a partition of $[0,M]$ into $k$ equal intervals but for  some  padding that we added to the leftmost and the rightmost intervals:
on the left we added $[-\frac{\log n}{n},0)$ and on the right $(M,M+1/n]$.
% $u^1_k:=$,   $u^i_k:=$ for $1<i<k$ and  $u^k_k:=$. 

For each $\mu\in\C$, $n,k\in\N$, $i=1..k$ define the sets 
\begin{equation}\label{eq:t}
T_{\mu,k,i}^n:=\left\{x_{1..n}\in \X^n: {1\over n}\log \frac{\mu(x_{1..n})}{\rho(x_{1..n})}\in u_k^i\right\}.
\end{equation} 
Observe that,  for every  $\mu,k,n\in\N$, these sets  constitute a partition of % $\X^n$
$T_{\mu}^n$ into $k$ disjoint sets: indeed, on the left we have ${1\over n}\log \frac{\mu(x_{1..n})}{\rho(x_{1..n})}\ge -{1\over n}\log n$ 
 by definition~\eqref{eq:tt}
of $T_{\mu}^n$, and on the right we have ${1\over n}\log \frac{\mu(x_{1..n})}{\rho(x_{1..n})}\le M+1/n$ from~\eqref{eq:boundedness}.
In particular, from this definition,    %Note that by construction~\eqref{eq:t} 
for all $x_{1..n}\in T^n_{\mu,k,i}$ we have
 \begin{equation}\label{eq:cops1}
 % 2^{{(i-1)M\over k}n}\rho(x_{1..n})\le 
 \mu(x_{1..n}) \le 2^{{iM\over k}n+1}\rho(x_{1..n}). %%% +1 from the case i=1 (say this?)
 \end{equation}
%{\em Step 2n: a countable cover, time $n$.}

%Fix an $n\in\N$.
 For every $n,k\in\N$ and $i\in\{1..k\}$ consider the following construction. 
 Define $m_{1}:=\max_{\mu\in\C}\rho(T_{\mu,k,i}^n)$ (since $\X^n$ are finite all suprema are reached). 
 Find any $\mu_{1}$ such that $\rho(T_{\mu_1,k,i}^n)=m_{1}$ and let
$T_{1}:=T^n_{\mu_1,k,i}$. For $l>1$, let $m_{l}:=\max_{\mu\in\C}\rho(T_{\mu,k,i}^n\backslash T_{l-1})$.
 If $m_l>0$, let $\mu_l$ be any $\mu\in\C$ such 
that $\rho(T^n_{\mu_l,k,i}\backslash T_{l-1})=m_l$, and let $T_l:=T_{l-1}\cup T^n_{\mu_l,k,i}$; otherwise let $T_l:=T_{l-1}$ and $\mu_l:=\mu_{l-1}$. 
%Note that, since $\mu\in\C$, we must have $T^n_{\mu,k,i}=\cup_{l\in\N}T^n_{\mu_l,k,i}$. Therefore, 
Note that,  for each $x_{1..n}\in T_{l}$ there is $l'\le l$ such that $x_{1..n}\in T^n_{\mu_{l'},k,i}$ and thus from~\eqref{eq:t} we get      %(similar to~\eqref{eq:cops1})
 \begin{equation}\label{eq:cops2}
  2^{{(i-1)M\over k}n-\log n}\rho(x_{1..n})\le \mu_{l'}(x_{1..n}).  %%%%\log n from the case i=1
   % \le 2^{{iM\over k}n}\rho(x_{1..n}).
 \end{equation}
 Finally, define 
\begin{equation}\label{eq:nun}
\nu_{n,k,i}:=\sum_{l=1}^{\infty} w_l\mu_l.
\end{equation}
(Notice that 
for every $n,k,i$ there is only a finite number of positive $m_l$,
since the set $\X^n$ is finite; thus the sum in the last definition is effectively finite.)
%  let $K_{n,k,i}$ be the largest index $l$ such that $m_l>0$.
% where by a slight abuse of notation we assume that the weights $w_l$ are renormalized to sum to~1.
We will show that the set $\{\nu_{n,k,i}:n,k\in\N, i=1..k\}$ is the countable set (sequence) which we are looking for to establish~\eqref{eq:nuro}.
Thus, we shall define the predictor $\nu$ as 
\begin{equation}\label{eq:nu}
\nu:=\sum_{n,k\in\N}w_nw_k{1\over k}\sum_{ i=1}^k\nu_{n,k,i},
\end{equation}
and show that~\eqref{eq:nuro} holds for $\nu$ so defined.

First we want to show that for each $\mu\in\C$ for each fixed $k,i$  the sets $T^n_{\mu,k,i}$ are covered up to a negligible  $\mu$-probability by the sets $T_l$ with indices  $l$ that are not too small. 
Observe that, by definition, for each $n,i,k$ the sets $T_l\backslash T_{l-1}$ are disjoint (for different $l$) and have non-increasing (with $l$) $\rho$-probability. Therefore, $\rho(T_{l+1}\backslash T_{l})\le 1/l$  for all $l\in\N$. Hence, from the definition of $T_l$, we must also have $\rho(T^n_{\mu,k,i}\backslash T_{l})\le 1/l$   for all $l\in\N$.
From the latter inequality and~\eqref{eq:cops1} we obtain $\mu(T^n_{\mu,k,i}\backslash T_{l})\le 2^{{iM\over k}n+1}\, 1/l$.
Consequently, for any $a>M/k$  taking $l:=2^{({iM\over k}+a)n+1}$ we obtain that for each $x_{1..n}\in T^n_{\mu,k,i}$  %%% typo k--l ?
 except possibly for a set of $\mu$-probability~$2^{-an}$ (that is, for $x_{1..n}\in T^n_{\mu,k,i}\backslash T_{l}$) there is $l'\le l$ such that  the following chain holds
\begin{equation}\label{eq:fnu}
% \nu(x_{1..n})\ge w_nw_kw_i w_{2^{({i\over k}+a)n}}\mu(x_{1..n})\ge w_nw_kw_i {2^{-({i\over k}+a)n}}\mu(x_{1..n}) \ge w_nw_kw_i {2^{-(a-{1\over k})n}}\rho(x_{1..n})
 \nu(x_{1..n})\ge w_nw_k{1\over k} w_{2^{({iM\over k}+a)n}}\mu_{l'}(x_{1..n})\ge 2^{-({iM\over k}+a)n+o(n)}\mu_{l'}(x_{1..n}) \ge 2^{-(a+{M\over k})n+o(n)}\rho(x_{1..n})
\end{equation}
where the first inequality is from~\eqref{eq:nu} and~\eqref{eq:nun} (with the value of $l$ we selected), the second is 
by definition of $w_l$ and the third uses~\eqref{eq:cops2}.

Suppose that there exist $\mu\in\C$ and $\delta'>0$ such that ${1\over n}d_n(\mu,\nu)> {1\over n}d_n(\mu,\rho) + \delta'$ infinitely often, so that 
${1\over n}\E_\mu\log\frac{\rho(x_{1..n})}{\nu(x_{1..n})}>\delta' $ i.o. Using~\eqref{eq:boundedness} %Since $-\log\nu(x_{1..n})\le n\log |\X|+const$ for all $x_{1..n}\in\X^n$ % ${1\over n}d_n(\mu,\nu)\le M+const$,
%from some $n$ on, 
we conclude that
%this implies that 
there exist $\epsilon',\delta>0$, an infinite sequence of indices $(n'_j)_{j\in\N}$ and   sets $A'_j\subset X^{n'_j}$ such that $\mu(A'_j)>\epsilon'$ and
 $-\log\nu(x_{1..n_j}) > {n_j\delta} -\log\rho(x_{1..n_j})$ for $x_{1..n_j}\in A'_j$.  % and let $\gamma=1/k$.
Taking into account~\eqref{eq:tm}, we also obtain $\mu(A'_j\cap T_{\mu}^n) > \epsilon'-1/n$.
Recall that for each $k\in\N$ the sets $T^n_{\mu,k,j}$ partition  each of the sets  $T_{\mu}^n$ and therefore 
% the sets $T^n_{\mu,k,j}$ also partition 
 each of the sets $A'_j\cap T_{\mu}^n$ into at most $k$ sets.
%Since the set $U_m$ is finite, 
Therefore, for every $k$ there must exist a cell of this partition, that is,  an index  $i\in\{1..k\}$, along with an $\epsilon\ge\epsilon'/k>0$
and  subsequences $(n_j)_{j\in\N}$ and $(A_j)_{j\in\N}$ (with  $A_j\subset X^{n_j}$)    of the sequences $(n'_j)_{j\in\N}$ and  $(A'_j)_{j\in\N}$     such that $\mu(A_j\cap T^{n_j}_{\mu,k,i})>\epsilon-1/n$ for 
all $j\in\N$.
Denote $B_j=A_j\cap T^{n_j}_{\mu,k,i}$ for each $j\in\N$. 
% and consider any set $B$ of this sequence, along with the corresponding index $n$ from the sequence $(n_i)_{i\in\N}$.
We have thus obtained, finally, an infinite sequence of indices $(n_j)_{j\in\N}$ and sets $B_j\subset  T^{n_j}_{\mu,k,i}$ of $\mu$-probability bounded from below 
by $\epsilon/2$, such that for each $x_{1..n_j}\in B_j$ we have 
\begin{equation}\label{eq:unnu}
  \nu(x_{1..n_j})< 2^{-\delta n_j} \rho (x_{1..n_j}).
\end{equation}
Take $k>0$ such that $M/k<\delta/4$. To conclude the proof of the second statement of the theorem, it remains to observe that~\eqref{eq:unnu} contradicts~\eqref{eq:fnu} with $a=\delta/2$.

Let  $\gamma_j>V_\C$, $j\in\N$ be a sequence such that $\lim_{j\to\infty}\gamma_j = V_\C$. 
Find then a sequence $\rho_j\in\mathcal P$ such that $D(\C,\rho_j)\le\gamma_j$.
Fix any $\rho\in\{\rho_j:j\in\N\}$.  We need to show that
\begin{equation}\label{eq:nuro2}
  D(\C,\nu)\le D(\C,\rho).
\end{equation}
So far  we have shown that for every $\rho_j$, $j\in\N$ there is  a measure $\nu_j$ 
of the form $\sum_{k\in\N}w_k'\mu_k$, where $\mu_k\in\C$
such that  $D(\C,\nu_j)\le D(\C,\rho_j)$. 
It remains to define $\phi:=\sum_{j\in\N} w_j\nu_j$ and show that it satisfies~\eqref{eq:nuro2}.
Indeed, for every $\mu\in\C$ and every $j\in\N$
\begin{equation*}
d_n(\mu,\phi)={1\over n}E_\mu\log\frac{\mu(x_{1..n})}{\phi(x_{1..n})}\le {1\over n}E_\mu\log\frac{\mu(x_{1..n})}{\nu_j(x_{1..n})}-{1\over n}\log w_j,
\end{equation*}
so that $D(\mu,\phi)\le D(\mu,\nu_j)\le D(\mu,\rho_j)\le\gamma_j$. Finally, recall that $\gamma_j\to V_\C$ to obtain the desired statement.
\end{proof}

%\section{Open questions}\label{s:disc}
\section{Decision-theoretic interpretations}\label{s:dec}
Classical decision theory is concerned with single-step games. Among its key  results 
are the complete class and minimax theorems. 
%The minimax and similar quantities arise in decision theory, which is concerned with single-step games.
The infinite-horizon case studied here presents both differences and similarities which we attempt  to summarize here. 
A distinction worth mentioning at this point is that the results presented here are obtained  under no assumptions whatsoever, whereas
the results in decision theory we refer to always have a number of conditions; on the other hand, here we are concerned with just one spcific
loss function (KL divergence) rather than general losses as is common in decision theory.
% The definition of the minimax value~\eqref{eq:vc} has a distinct decision-theoretic flavour to it, and in particular one can try to ask 
% whether the minimax  and the complete class theorems hold for it. Here we make the flavour and the question precise. 

Predictors $\rho\in\mathcal P$ are called  {\em  strategies  of the statistician}. 
The measures $\mu\in\C$ are now the basic {\em  strategies of the opponent}, and the first thing  we need to do is to extend these to  randomized strategies. 
To this end, denote $\C^*$ the set of all probability distributions over measurable subsets of $\C$. Thus, the opponent selects a randomized
strategy $W\in\C^*$ and the statistician (predictor) $\rho$ suffers the loss 
\begin{equation}\label{eq:ew}
  E_{W(\mu)} D(\mu,\rho),
\end{equation}
 where the notation $W(\mu)$ means that $\mu$ is drawn
according to $W$. Note a distinction with the combinations we considered before. A combination of the kind $\nu=\int_\C dW$ is itself
a probability measure over the one-way infinite sequences, whereas a measure $W\in \C^*$ is a measure over $\C$.
In other words, the difference is between putting the integral $\int_\C dW$ outside of $D$ as in \eqref{eq:ew} or inside of $D$
which would be $D(\int_\C dW(\mu),\rho)$. In the terminology of \cite{Gray:88},
the measure $\int_\C dW(\mu)\in\mathcal P$ is the barycentre of $W\in\C^*$.

{\noindent\bf Minimax.} Generalizing the definition~\eqref{eq:vc} of $V_\C$, we can now introduce the {\em upper value}
\begin{equation}\label{eq:bvc}
   \bar V_\C:=\inf_{\rho\in\mathcal P} \sup_{\mu\in\C^*} E_{W(\mu)} \bd(\mu,\rho).
\end{equation}
Furthermore, the {\em maximin} (the {\em lower value}) is defined as 
\begin{equation}\label{eq:ubvc}
   {\underline V}_\C:= \sup_{\mu\in\C^*} \inf_{\rho\in\mathcal P} E_{W(\mu)} \bd(\mu,\rho).
\end{equation}

The so-called minimax theorems in decision theory (e.g., \cite{Ferguson:14}) for single-step games and general loss functions state that,
 under certain conditions, $\bar V_\C=\underline V_\C$ and the 
statistician has a minimax strategy, that is, there exists $\rho$ on which $\bar V_\C$ is attained. 
Minimax theorems generalize the classical result of von~Neumann~\cite{neumann:28}, and provide sufficient conditions of various generality for it to hold.
A rather  general sufficient condition is the existence of a topology with respect to which the set of all strategies of the statistician, $\mathcal P$  in our case, 
is compact, and the risk, which  is $D(\mu,\rho)$ in our case, is lower semicontinuous. 
Such a condition seems nontrivial to verify. For example, a (meaningful) topology with respect to 
which $\mathcal P$ is compact is that of the so-called distributional distance \cite{Gray:88} (in our case it coincides with the topology of the 
weak${}^*$ convergence), but $D(\mu,\rho)$ is not  (lower) semicontinuous with respect to it. 
Some other (including non-topological) sufficient conditions  are given in \cite{Sion:58,Lecam:55}.

In our setup, it is easy to see that $\bar V_\C = V_\C$ and so Theorem~\ref{th:2} holds for $\bar V_\C$. 
Thus, using decision-theoretic terminology, we can state the following.
\begin{corollary}[minimax]
 For every set $\C$ of strategies of the opponent, the statistician has a minimax strategy.
\end{corollary}

% In our setup, from Theorem~\ref{th:2} we know that the statistician has a minimax strategy. So 

However, the question  of whether
 the upper and the lower values coincide remains open. 
That is, we are taking the worst possible distribution over $\C$, and ask what is the best possible 
predictor  with the knowledge of  this distribution ahead of time. The question is whether ${\underline V}_\C = V_\C$.
A closely related question is  whether  there is a worst possible strategy for the opponent.
 This latter would be somehow a maximally spread-out (or maximal entropy) distribution over $\C$.
In general, measurability issues  seem to be very relevant here, especially for the maximal-entropy distribution part.

% , and I would not conjecture, without any additional assumptions, 
% that the minimax theorem holds. 
{\noindent\bf Complete class.}
For a set of measures (strategies of the opponent) $\C$, a  predictor  %(rule/ strategy of the statistician) 
$\rho_1$ is said to be  {\em as good as}  a predictor $\rho_2$
if $D(\mu,\rho_1)\le D(\mu,\rho_2)$ for all $\mu\in\C$. A predictor $\rho_1$ is {\em better (dominates)} $\rho_2$ if $\rho_1$ is as good as $\rho_2$
and $D(\mu,\rho_1) < D(\mu,\rho_2)$ for some $\mu\in\C$. 
A predictor $\rho$ is {\em  admissible} (also called {\em Pareto optimal}) if there is no predictor $\rho'$ which is better than $\rho$; 
otherwise it is called {\em inadmissible}.
Similarly, a set of predictors $D$ is called a {\em complete class} if for every $\rho'\notin D$ there  is $\rho\in D$ such that 
$\rho$ is better than $\rho'$. A set of of predictors $D$ is called an {\em  essentially complete class} if 
for every $\rho'\notin D$ there  is $\rho\in D$ such that 
$\rho$ is as good as $\rho'$. 
An (essentially) complete class is called {\em minimal} if none of its proper subsets is (essentially) complete.

Furthermore, in  decision-theoretic terminology, a predictor $\rho$ is called a {\em Bayes rule} for a prior $W\in\C^*$ 
if it is optimal for $W$, that is, if it attains $\inf_{\rho\in\mathcal P} E_{W(\mu)}D(\mu,\rho)$. 
Clearly, if $W$ is concentrated on a finite or countable set then any mixture over this set with full support
is a Bayes rule, and the value of the $\inf$ above is 0.

In decision theory, the complete class theorem (\cite{Wald:50,Lecam:55}, see also  \cite{Ferguson:14}) states that, under certain conditions similar to those
above for the minimax theorem, the set of Bayes rules is complete and the admissible Bayes rules  form a minimal complete class.
%  This means that 
% there is no predictor that is better than some Bayes rule on all $W\in\C^*$. 

An important difference in our set-up is that all strategies are inadmissible (unless $V_\C$=0), and one cannot speak about 
minimal (essentially) complete classes. However, the set of all Bayes rules 
is still essentially complete, and an even stronger statement holds: it is enough to consider all Bayes rules with countable priors:

\begin{proposition}
 For every set $\C$, the  set  of those Bayes rules %(for all $W\in\C^*$) 
 whose priors are concentrated on at most countable sets is essentially complete.
 There is no admissible rule (predictor) and no minimal essentially complete class
 unless $V_\C=0$. In the latter case, every 
predictor $\rho$ that attains this value is admissible and the set $\{\rho\}$ is minimal essentially complete. 
 %There is no minimal essentially complete class unless $V_\C=0$.
\end{proposition}
\begin{proof}
The first statement is a reformulation of the second statement of Theorem~\ref{th:2}.
 To prove the second statement, consider any $\C$ such that $V_\C>0$, take a predictor $\rho$ that attains this value (such a predictor exists by Theorem~\ref{th:2}),
and a measure $\mu$ such that $D(\mu,\rho)>0$. Then for a predictor $\rho':=1/2(\rho+\mu)$ we have $D(\mu,\rho')=0$, so that $\rho'$ is better than $\rho$
and thus $\rho$ is inadmissible.  The statement about minimal essentially complete class is proven analogously.
The statement about the case $V_\C=0$ is obvious.
\end{proof}

% 
% 
% Theorem~\ref{th:2} above then establishes  the complete class theorem in a somewhat stronger form: 
% it says that the class of all Bayes rules for all $W\in\C^*$ {\em whose priors are concentrated on at most countable
% sets} is essentially complete. (In fact, one of these Bayes rules is minimax.)
% 
% Thus, for our prediction problem, we have, in Theorem~\ref{th:2}, the complete class theorem in a somewhat stronger form, and 
% about a half of the minimax theorem: the statistician has a minimax strategy, but we don't know whether 
% $\bar V_\C=\underbar V_\C$. 

\section{Examples}\label{s:ex}
In \cite{Ryabko:10pq3+} several examples are considered in detail for the case $V_\C=0$; these include 
the case of countable $\C$, the set of i.i.d.\ measures, Markov chains, bounded-memory processes and stationary ergodic 
processes. Therefore, here we will only look at the case $V_\C>0$. For simplicity, we assume $\X=\{0,1\}$ in the examples.

\noindent{\bf Typical Bernoulli 1/3 sequences.}
We start with a somewhat artificial example, but a one on which it is relatively easy to see how countable mixtures give predictors
for large uncountable sets.
%Consider the following example.
 Take the binary $\X$ and consider all sequences $\x\in\X^\infty$ such 
that the limiting number of 1s in $\x$ equals $1/3$. Denote the set of these sequences $S$ and let the set $\C$ consist of all Dirac measures concentrated on sequences from $S$.   Observe that the Bernoulli i.i.d.\ measure $\delta_{1/3}$ with 
probability $1/3$ of 1 predicts measures in $\C$ relatively well:  $D(\C,\delta_{1/3})=h(1/3)$,
where $h$ stands for the binary entropy, and this is also the minimax loss for this set, $V_\C$. 
It might then appear surprising that this loss is achievable by a combination of  countably many measures from $\C$~--- after all,  this
set consists only of deterministic measures. Let us try to see what such a combination may look like. 
By definition, for any sequence $\x\in S$ and every $\epsilon$ we can find $n_\epsilon(\x)\in\N$ such that for all $n\ge n_\epsilon(\x)$
the average number of 1s in $x_{1..n}$ is within $\epsilon$ of $1/3$.  Fix the sequence of indices $k_j:=2^j$, $j\in\N$ and the sequence
of thresholds   $\epsilon_l:=2^{-l}$.  For each $k_j$ let ${S'}_j^l\subset S$  be the  set of all sequences $\x\in S$ such that $n_{\epsilon_l}(\x)<n_j$. Select then a finite subset $S_j^l$ of ${S'}_j^l$ such that for each $\x'\in {S'}_j^l$ there is $\x\in S$ such 
that $x'_{1..n_j}=x_{1..n_j}$. This is of course possible since the set $\X^{n_j}$ is finite. Now for each $\x\in S_j^l$ take 
the corresponding measure $\mu_\x\in\C$ and attach to it the weight $w_lw_j/|S_j^l|$, where, as before, we are using 
the weights $w_k=w/k\log^2k$. Taking these measures for all $j,l\in\N$ we 
obtain our convex combination.  Of course we did not enumerate all sequences in $S$ (or measures in $\C$) this way; but for 
each sequence $\x\in S$ and for each $n$ there is a sequence among those that we did enumerate that coincides with $\x$ up to the index $n$.   One can then use the theory of types \cite{Csiszar:98} to calculate the sizes of the sets $S_j^l$
and to check that the weights we found give the optimal loss we are after; but for the illustrative purposes of this example 
this is already not necessary.

\noindent{\bf Processes with  abrupt changes.}
Start with a family of  distributions $S$, for which we have a good predictor, for example $S$ is the set $B$ of all Bernoulli i.i.d.\ processes, or 
more generally a set for which $V_S=0$.  The family $\C_\alpha$ parametrized by $\alpha\in(0,1)$ and $S$  is then  the family of all process constructed as follows:  there is a sequence of indexes $n_i$ such that $X_{n_i..n_{i+1}}$ is distributed according to $\mu_i$ for some $\mu_i\in S$. Take  then all possible sequences $\mu_i$ and all sequences $n_i$ whose limiting  frequency $\lim_{i\to\infty}\{i:n_i<n\}$ is bounded by $\alpha$  to obtain our set $\C_{S,\alpha}$. Thus, we have a family of processes with abrupt changes in 
distribution, where between changes the distribution is from $S$, the changes are assumed to have the frequency bounded by $\alpha$ but are otherwise arbitrary.
This example was considered in \cite{Willems:96} for the case $S=B$, with the goal of minimizing the regret w.r.t.\ the predictor that knows where the changes
occur (the value $V_\C$ was not considered directly). The  method proposed in the latter work, in fact, is not limited to the case $S=B$, but is general. The algorithm is based on a prior over all possible sequences $n_i$ of changes; between the changes
the optimal predictor for $B$ is used, which is also a Bayesian predictor with a specific prior. The regret obtained is of order $\log n$. Since for Bernoulli processes themselves  the best achievable loss up to time $n$ is  $1/2\log n +1$, we can see that $V_{\C_{B,\alpha}}=\alpha(1-1/2\log\alpha)$.
A similar result can be obtained if we replace   Bernoulli processes with Markov processes, but not with an arbitrary $S$ for which $V_S=0$. For example, if we take $S$ to be all finite-memory 
distributions, then the resulting process may be completely unpredictable ($V_\C=1$): indeed, if the memory of distributions $\mu_i$ grows (with $i$) faster than $\alpha n$,
then there is little one can do.  For such sets $S$ one can make the problem amenable by restricting the way the distributions $\mu_i$ are selected, 
for example, imposing an ergodicity-like condition that the average distribution has a limit. Another way (often considered in the literature in 
slightly different settings, see \cite{Gyorgy:12} and references therein) is to have $\alpha\to0$, although in this case one recovers $V_{\C_S}=0$ provided 
$\alpha$ goes to 0 slowly enough (and, of course, provided $V_S=0$).

\noindent{\bf Predictable aspects.} 
% We continue with examples of situations where asymptotically accurate prediction is not possible: the best achievable asymptotic average error 
% is positive ($V_\C>0$), but still one can make a lot of sense out of the data and learn to predict quite well ($V_\C<\log|X|$).
% The common conclusion is that in the examples described below it is possible to construct a minimax predictor by a (countable) Bayesian mixture. However, we will not present the details of such predictors, and just consider examples (or example classes) of the sets $\C$ instead.
The preceding example can be thought of as an instantiation of the general class of processes in which  some aspects are predictable while 
others are not. Thus, in the considered example changes between the distributions were unpredictable, but between the changes the distributions were predictable. % A rich class of examples comes from the  idea that some aspects of the process may be predictable while others are not.
% The simplest toy example\footnote{This example was used in \cite{Hutter:04uaibook} for a different purpose.} is that of all distributions such that on  even steps the outcomes are Bernoulli i.i.d.\ but 
% on odd steps are arbitrary (that is, may form an arbitrary deterministic sequence). 
%  If $\C$ is the class of all such processes then $V_\C=1/2$.  One can replace Bernoulli i.i.d.\ here with  Markov chains, or any other set of 
% processes  $B$ for which $V_B=0$. 
%
%A more meaningful generalization is as follows. 
Another example of this kind is that of processes predictable on some scales but not on others. 
 Imagine that it is  possible to predict, for example, large fluctuations of the process but not small fluctuations (or the other way around). More formally, consider now an alphabet $\X$ with $|\X|>2$, and
let $Y$ be some partition of $\X$. For any sequence ${x_1,\dots,x_n,\dots}$ there is an associated sequence $y_1,\dots,y_n,\dots$ where 
$y_i$ is defined as $y\in Y$ such that $x_i\in y$.  Here again we can obtain examples of sets $\C$ of processes with $V_\C\in(0,1)$ by restricting the distribution of $y_1,\dots,y_n,\dots$ to a set $B$ with $V_B=0$. The interpretation is that, again, we can model 
the $y$ part (by processes in $B$) but not the rest, which we then allow to be arbitrary.

Yet another example is that of  processes  predictable only after certain kind of events: for example, after a price drop; or after a rain. At other times, the process is unpredictable: it can, again, be an arbitrary deterministic sequence.   More formally, let a set $A\subset \X^*:=\cup_{k\in\N}\X^k$ be measurable. Consider for each sequence  $\x={x_1,\dots,x_n,\dots}$
  another (possibly finite) sequence $\x'={x'_1,\dots,x'_n,\dots}$ given by $x'_i:=(x_{n_i+1})_{i\in\N}$ where $n_i$ are all indexes such that $x_{1..n_i}\in A$. 
We now form the set $\C$ as the set of all processes $\mu$ such that $\x'$ belongs ($\mu$-a.s.) to some pre-defined set $B$;
 for this set $B$ we may have  $V_B=0$.  This means that we can model what happens after events in $A$~--- by processes in $B$, but not the rest of the times, on which we say the process may be arbitrary. For different $A$ and $B$ we then obtain  examples where $V_\C\in(0,1)$.

 %Another way to look at predictable aspects is as follows.
\section{Relation to the non-realizable case}\label{s:not}
As mentioned in the Introduction, \cite{Ryabko:16pqnot} shows that in the non-realizable case all Bayesian mixture predictors
may be suboptimal. Here we make this statement precise in order to clarify its relation to the main result. 

The non-realizable case is when the measure generating the data does not belong to $\C$. We are then looking 
at the set $\C$ as the set of experts or models, and we seek a predictor $\rho$ that predicts any measure $\nu$ (that generates the data) whatsoever  as well as 
the best (for this $\nu$)  $\mu\in\C$. 

Thus, if we have two predictors $\mu$ and $\rho$, we can  define the {\em regret} up to time $n$ of 
(using the predictor) $\rho$ as opposed to (using the predictor) $\mu$ on the measure $\nu$ (that is, $\nu$ generates the sequence to predict) as 
$$
R_n^\nu(\mu,\rho):=d_n(\nu,\rho) - d_n(\nu,\mu).
$$

Furthermore, define the asymptotic average regret as 
$$
\bar R^\nu(\mu,\rho):=\limsup_{n\to\infty}{1\over n}R^\nu_n(\mu,\rho),
$$
and
$$
\bar R^\nu(C,\rho):= \sup_{\mu\in C}\bar R^\nu(\mu,\rho).
$$

It is shown in~\cite{Ryabko:16pqnot} that there exists a set $\C$ such that any  Bayesian predictor must have a linear regret, while there exists a predictor with a sublinear regret:
\begin{theorem*}[\cite{Ryabko:16pqnot}]
  There exist a set $C$ of measures and a predictor $\rho$ such that 
%$R_T^\nu(C,\rho)\le cT +o(T)$, 
% $\bar R^\nu(C,\rho)\le c$ for every measure $\nu$, yet
 $\bar R^\nu(C,\rho) = 0 $ for every measure $\nu$, yet
for every Bayesian predictor $\phi$ with a prior concentrated on $C$ 
there exists a measure $\nu$ such that  % $R_T^\nu(C,\phi)\ge 2cT + o(T)$ and 
we have $\bar R^\nu(C,\phi)\ge c>0$ 
where $c$ is a  constant (independent of $\phi$).
%In other words, any  Bayesian predictor must have a linear regret, while there exists a predictor with a sublinear regret.
\end{theorem*}
\cite{Ryabko:16pqnot} also argues that this applies more broadly than just Bayesian predictors: all meaningful combinations of measures 
in~$\C$ may be useless for minimizing regret. 
We remind again that such a set $\C$ must necessarily be uncountable.

% We thus reach the conclusion that it is better to make the model class large enough to be sure it includes the ``true'' measure  generating the data 
% and seek to minimize the loss for this (large) set, rather than minimizing regret for a smaller set that may not include the true measure.
% 

\section{Discussion}
A statistician facing an unknown stochastic phenomenon has a large, nonparametric model class at hand that 
she has reasons to believe captures some aspects of the problem. Yet other aspects remain completely enigmatic, 
and there is little hope that the process generating the data indeed comes from the model class.
For this reason the statistician is content at having non-zero error no matter how much data may become available 
now or in the future, but she would still like to make some use of the model. There are now two rather distinct ways 
to proceed. One  is to say that the data may come from an arbitrary deterministic sequence, and try to construct 
a predictor that minimizes the regret with respect to every distribution in the model class, on every deterministic sequence. The other way is to try to enlarge
the model class, in particular, by allowing that all there is enigmatic in the process may be arbitrary (that is, 
an arbitrary deterministic sequence). This second  way may be more difficult precisely on the modelling step. Yet, 
the conclusion of this work is that this is the way to follow, for in this case one can be sure that it is possible
to make statistical inference by standard available tools, specifically, Bayesian forecasting:  even if 
the best achievable asymptotic error is non-zero it is attained by a Bayesian forecaster with some prior.
Finding such a prior is a separate problem, but it is a one with  which Bayesians are familiar. Here, modelling that ``enigmatic''
part should not create much trouble:  a good distribution  over all deterministic sequences is just the Bernoulli i.i.d.\ measure
with equiprobable outcomes.
(Note that it is not necessary to look for priors concentrated on  countable sets.) 
On the other hand, for the regret-minimization route, the statistician cannot use an arbitrary model class; indeed, 
she would first need to make sure that regret minimization is viable  at all for the model class at hand: it may happen that every combination 
of distributions in the model is suboptimal. There are no criteria for checking this, only some (rather small) examples, such as finite or countable sets, or specific parametric families.

Finding such criteria for the viability of  regret minimization is  an interesting open problem.
To make it more precise, the question is for which sets $\C$ of distributions the minimax regret 
(is attainable and) can be attained by a combination (either Bayesian or some other) of distributions in~$\C$.

It is worth noting that the conclusions of the paper are not about Bayesian versus non-Bayesian inference. 
Rather, Bayesian inference is used as a generic approach to construct predictors for general (uncountable) model
classes.  At this level of generality it is hard to find any  alternative  approach, although it would 
be interesting to see which predictors can be generalized (to arbitrary model classes) and whether the corresponding
result holds for them. The negative result of \cite{Ryabko:16pqnot}, as explained there, is not restricted
to Bayesian predictors but holds in any foreseeable  generality.

Another interesting open  question concerns different losses. While the proof does not seem to be hinged very specifically 
on the log loss, it does use some properties of it in an important way. In particular, the property that if $\mu$ predicts
$\nu$ then also any convex combination $\alpha\mu+(1-\alpha)\rho$ predicts  $\nu$ for any $\rho$. This does not hold
for some other losses, in particular already for KL loss without Cesaro averaging; see \cite{Ryabko:08pqaml} for a discussion
and some results on this property.  

Some other interesting open-question are the decision-theoretic ones mentioned in Section~\ref{s:dec}; specifically, those
concerning the minimax theorem and the existence of  maximally spread distributions over $\C$.

Finally, an intriguing question is whether a result like Theorem~\ref{th:2} holds if one allows convergence rates into consideration. Now 
that we know that the minimax asymptotic error is achievable, we can ask whether the minimax rate of convergence 
to this error is also achievable (by a Bayesian predictor). The proof of the version of Theorem~\ref{th:2} for the $V_\C=0$ case in \cite{Ryabko:10pq3+} 
clearly does not generalize to achieve such a result (the rates one extract from that proof are rather bad), but with the present proof this may be possible. 
%\subsection{Other losses} 
%\subsection{}

%\newpage
\bibliographystyle{plain}

\end{document}